\documentclass{amsart}

\usepackage{amsmath,amstext,amsfonts}

\begin{document}
\title[Singly generated composition $C^*$-algebras]{Composition operators within singly generated composition $C^*$-algebras}

\author{Thomas L. Kriete}
\address{Department of Mathematics,
University of Virginia, Charlottesville, VA 22904}
\email{tlk8q@virginia.edu}

\author{Barbara D. MacCluer}
\address{Department of Mathematics,
University of Virginia, Charlottesville, VA   22904}
\email{bdm3f@virginia.edu}

\author{Jennifer L. Moorhouse}
\address{Department of Mathematics,
Colgate University, Hamilton, NY   13356}
\email{jmoorhouse@mail.colgate.edu}


\date{July 10, 2008}

\newtheorem{thm}{Theorem}
\newtheorem{prop}{Proposition}
\newtheorem{lemma}{Lemma}
\newtheorem{cor}{Corollary}

\newcommand{\p}{\varphi}
\newcommand{\cp}{C_{\varphi}}
\newcommand{\D}{{\mathbb D}}
\newcommand{\B}{{\mathcal B}}
\newcommand{\com}{{[C_{\varphi}^*,C_{\varphi}]}}

 \newcommand{\ebox}{\hspace*{1em}\linebreak[0]\hspace*{1em}\hfill\rule[-1ex]{
.1in}{.1in}\\}
 \newcommand{\pf}{\noindent \bf Proof.\rm\ \ }

\begin{abstract}
Let $\varphi$ be a linear-fractional self-map of the open unit disk $\D$, not an automorphism,
such that $\varphi(\zeta)=\eta$ for two distinct points $\zeta,\eta$ in the unit circle
$\partial \D$.  We consider the question of which composition operators lie
in $C^*(\cp,{\mathcal K})$, the unital $C^*$-algebra generated by the composition operator
$C_{\varphi}$ and the ideal ${\mathcal K}$ of compact operators, acting on the Hardy
space $H^2$.  This necessitates a companion study of the unital $C^*$-algebra generated by
the composition operators induced by all parabolic non-automorphisms with common
fixed point on the unit circle.
\end{abstract}

\maketitle

\section{Introduction}

Given any analytic self-map $\varphi$ of the unit disk $\D$ in the complex plane, one can 
form the composition operator $C_{\varphi}:f\rightarrow f\circ\varphi$, which acts as a
bounded operator on the Hardy space $H^2$.  This paper is the second in a series
of three investigating spectral theory in $C^*$-algebras generated by certain composition 
and Toeplitz operators. In the first article \cite{kmm}, we
studied
$C^*(T_z,C_{\varphi})$, the unital $C^*$-algebra generated by the unilateral shift
$T_z$ on $H^2$ and a single composition operator $C_{\varphi}$ with $\varphi$ satisfying
\begin{equation}\label{defofph}
\left\{\begin{array}{l}
\varphi\mbox{ is a linear-fractional self-map of }\D\mbox{ which is not an automorphism, and}\\
\varphi(\zeta)=\eta\mbox{ for distinct points }\zeta, \eta\mbox{ in the unit circle }\partial\D.
\end{array}\right.
\end{equation}
Throughout the current paper, $\varphi$ will always have this meaning. 
The algebra $C^*(T_z,\cp)$ necessarily
contains the ideal ${\mathcal K}$ of compact operators on $H^2$.  The main result of
\cite{kmm} identifies $C^*(T_z,C_{\varphi})/{\mathcal K}$ with a certain $C^*$-algebra of $2\times 2$ matrix
valued functions; see Theorem 4.12 of \cite{kmm}.  The case where $\varphi$ is replaced 
by an automorphism
of $\D$, or even a discrete group of automorphisms, has been studied by Jury \cite{j}, \cite{j2}, 
and has a rather different character.

The shift $T_z$ does not appear to play a role in the questions we consider in this paper;
accordingly we omit it and study $C^*(C_{\varphi},{\mathcal K})$, the unital $C^*$-algebra generated
by $C_{\varphi}$, for $\varphi$ as described above, and the compact operators.  The composition
$\varphi\circ\varphi$ has sup-norm strictly less than $1$, so that
$C_{\varphi}^2=C_{\varphi\circ\varphi}$ is compact and non-zero.  Since
$\varphi$ has no boundary fixed point, a theorem of Guyker \cite{g}
shows that $C_{\varphi}$ is irreducible if and only if $\varphi(0)\neq 0$.
It follows that when $\varphi(0)\neq 0$
the unital $C^*$-algebra $C^*(C_{\varphi})$ generated by $C_{\varphi}$ alone contains ${\mathcal K}$;
see \cite{cocot}, p.74.
We want our $C^*$-algebras to always contain ${\mathcal K}$,  and we indicate this by continuing to write
$C^*(C_{\varphi},{\mathcal K})$ if the irreducibility criterion $\varphi(0)\neq 0$ is in
doubt.

Let ${\mathcal P}$ denote the dense subalgebra of $C^*(C_{\varphi},{\mathcal K})$
consisting of all finite linear combinations of the identity $I$, all words in
$C_{\varphi}$ and $C_{\varphi}^*$, and all compact operators.  An element $B$ in
${\mathcal P}$ has a unique representation of the form
\begin{equation}\label{defofB}
B=cI+f(\cp^*\cp)+g(\cp\cp^*)+\cp p(\cp^*\cp)+\cp^*q(\cp\cp^*)+K
\end{equation}
where $f,g,p$ and $q$ are polynomials, $f(0)=0=g(0)$, $c$ is complex, and $K$
is compact.
Let $s=1/|\varphi'(\zeta)|$ and write ${\mathcal D}$ for the $C^*$-algebra
of continuous $2\times 2$ matrix-valued functions $F$ on $[0,s]$ with
$F(0)$ a scalar multiple of the identity, equipped with the supremum operator
norm.  It was shown in \cite{kmm} that
there is a unique $^*$-homomorphism $\Psi$ of $C^*(\cp,{\mathcal K})$ onto ${\mathcal D}$
with $\mbox{Ker}\ \Psi={\mathcal K}$ and such that
\begin{equation}\label{defofpsi}
\Psi(B)=\left[\begin{array}{lr}
c+g&rp\\rq&c+f\end{array}\right]
\end{equation}
where $B$ is given by Equation~(\ref{defofB}) and $r(t)=\sqrt{t}$. Equivalently,
we have a short exact sequence of $C^*$-algebras
$$0\rightarrow \mathcal{K}\stackrel{i}{\rightarrow} C^*(\cp,{\mathcal K})\stackrel{\Psi}
{\rightarrow}{\mathcal D}\rightarrow 0$$
where $i$ is inclusion. For any operator $T$ on $H^2$ we write $\|T\|_e$ for
the essential norm of $T$, that is, the distance from $T$ to the ideal
${\mathcal K}$.
We note that if $T$ is in $C^*(\cp,{\mathcal K})$, then
$\|T\|_e=\|\Psi(T)\|$.

For bounded operators $A$ and $B$ on $H^2$, let us write
$A\equiv B\ (\mbox{mod }{\mathcal K})$ if there exists a compact operator
$K$ with $A=B+K$.
In \cite{kmm} the authors used C. Cowen's well-known adjoint formula \cite{cow2} to
show
that if 
$$\psi(z)=\frac{az+b}{cz+d}$$ is a linear-fractional self map of $\D$, not
an automorphism but satisfying $|\psi(z_0)|=1$ for some $z_0\in\partial\D$,
then
\begin{equation}\label{adjformula}
C_{\psi}^*\equiv\frac{1}{|\psi'(z_0)|}C_{\sigma_{\psi}}\ (\mbox{mod }{\mathcal K})
\end{equation}
where $\sigma_{\psi}$ is the so-called ``Krein adjoint" of $\psi$, 
$$\sigma_{\psi}(z)=\frac{\overline{a}z-\overline{c}}{-\overline{b}z+\overline{d}}.$$
We denote the Krein adjoint of $\varphi$ itself simply by $\sigma$, which the reader can distinguish
by context from the spectrum $\sigma(T)$ or essential spectrum $\sigma_e(T)$ of an operator
$T$.  

Now consider an operator $B$ in the dense subalgebra ${\mathcal P}$, as expressed in
Equation~(\ref{defofB}).  Since $I=C_z, C_{\varphi}^*C_{\varphi}\equiv sC_{\varphi\circ\sigma}\ 
(\mbox{mod }{\mathcal K})$ and $\cp C_{\varphi}^*\equiv s C_{\sigma\circ\varphi}\ (\mbox{mod } {\mathcal K})$,
where $s=1/|\varphi'(\zeta)|$, we can rewrite Equation~(\ref{defofB}) as 
\begin{equation}\label{newB}
B=cC_z+A_1+A_2+A_3+A_4+K'
\end{equation}
where 
$K'$ is compact and $A_1,A_2,A_3,$ and $A_4$ are finite linear combinations
of composition operators whose associated self-maps are chosen, respectively,
from the four lists
\begin{equation}\label{lists}
\{(\varphi\circ\sigma)_{n_1}\},\{(\sigma\circ\varphi)_{n_2}\},\{(\varphi\circ\sigma)_{n_3}
\circ\varphi\},\{(\sigma\circ\varphi)_{n_4}\circ\sigma\}.
\end{equation}
Here we write $(\psi)_n$ for the $n^{th}$ iterate of a self-map $\psi$ of $\D$, and let
$n_k$ range over the positive integers for $k=1,2$, and over the  non-negative integers for $k=3,4$.
See \cite{kmm} for further details.
Thus $C^*(\cp,{\mathcal K})$ is spanned, modulo the compacts, by actual composition
operators.  This leads to our main question:  for which analytic self-maps $\psi$ of $\D$ does
$C_{\psi}$ lie in $C^*(\cp,{\mathcal K})$?

In particular we will describe explicitly, in both function-theoretic
and operator-theoretic terms, all linear-fractional composition operators 
lying in $C^*(C_{\varphi},{\mathcal K})$.  This description plays a role in the third paper of
our series \cite{kmm3} which is devoted to spectral theory in Toeplitz-composition algebras
with several generators.  Along the way here we show that if $C^*({\mathbb P}_{\gamma})$
denotes the unital $C^*$-algebra generated by composition operators induced by the parabolic
non-automorphisms  fixing $\gamma$ in the unit circle, then there is a short exact sequence
$$0\rightarrow \mathcal{K}{\rightarrow} C^*({\mathbb P}_{\gamma})
{\rightarrow}C([0,1])\rightarrow 0$$
of $C^*$-algebras, where $C([0,1])$ denotes the algebra of continuous functions
on the unit interval.

\section{Necessary conditions}

If $\psi:\D\rightarrow \D$ is analytic we write $F(\psi)$ for the set of points
$\alpha$ in $\partial\D$ at which $\psi$ has a finite angular derivative in the
sense of Caratheodory; see \cite{cmbook},\cite{sbook}.  In particular, if $\alpha$ is in $F(\psi)$, the
nontangential limit $\psi(\alpha)$ necessarily exists and has modulus one. We write
$\psi'(\alpha)$ for the angular derivative of $\psi$ at $\alpha$.  Recall that this is
the ordinary derivative if $\psi$ extends analytically to a neighborhood of $\alpha$
and $|\psi(\alpha)|=1$.

It is well known that if $C_{\psi}$ is compact on $H^2$, then $F(\psi)$ is empty \cite{st}.
When $C_{\psi}$ is considered as acting on the Bergman space in $\D$, the converse assertion
is true \cite{ms}.  On the space $H^2$ considered here, however, ``$C_{\psi}$ is compact" is a strictly
stronger requirement than ``$F(\psi)$ is empty"; see, for example, \cite{cmbook},\cite{sbook}.  Our first goal
is to show that when $C_{\psi}$ lies in $C^*(\cp,{\mathcal K})$ these two conditions 
are equivalent. 

First we recall that a linear-fractional self-map $\rho$ of $\D$ is {\it parabolic} if $\rho$ fixes one point
$\gamma$ in $\partial\D$ and is conjugate, via the map $(\gamma+z)/(\gamma-z)$, to 
translation in the right half-plane $\Omega=\{w:\mbox{Re }w>0\}$ by a complex number $a$ with
non-negative real part. We denote this parabolic map by $\rho_a$, or by $\rho_{\gamma,a}$ if the fixed
point $\gamma$ is not clear from the context. A linear-fractional map $\rho$ 
with fixed point $\gamma$ in $\partial \D$ is parabolic provided $\rho'(\gamma)=1$.
Another representation of $\rho_a$ will prove useful.  The map $\tau_{\gamma}(z)=i(\gamma-z)/(\gamma+z)$
carries $\D$ onto the upper half-plane $\{w:\mbox{Im }w>0\}$ and takes $\gamma$ to
$0$, rather than infinity.  We write $u$ for the conjugate of $\rho_a$ by $\tau_{\gamma}$:
$u=\tau_{\gamma}\circ\rho_a\circ\tau^{-1}_{\gamma}$.  One readily computes that $u(w)=iw/(i+wa)$, and so
\begin{equation}\label{u}
u''(0)=2ia.
\end{equation}

Also important for us will be several lower bounds for the essential norm of a linear combination
of composition operators.  Given an analytic self-map $\psi$ of $\D$ and $\alpha$ in $F(\psi)$,
we call $D_1(\psi,\alpha)\equiv(\psi(\alpha),\psi'(\alpha))$ the {\it first order data vector} for
$\psi$ at $\alpha$.  If we have a finite collection
of maps $\psi_1,\psi_2\cdots,\psi_n$ and $\alpha$ lies in the union of the finite angular derivative sets
$F(\psi_1),F(\psi_2),\cdots F(\psi_n)$, we define
$${\mathcal D}_1(\alpha)=\{D_1(\psi_j,\alpha):1\leq j\leq n\mbox{ and } \alpha\in F(\psi_j)\},$$
the set of possible first order data vectors at $\alpha$ coming from that collection of maps.
Theorem 5.2 in \cite{km} states that if ${\mathcal D}_1(\alpha)$ is non-empty, then
for any complex numbers $c_1,\cdots, c_n$,

\begin{equation}\label{lb1}
\|c_1C_{\psi_1}+\cdots+c_nC_{\psi_n}\|_e^2\geq \sum_{{\bf d}\in{\mathcal D}_1(\alpha)}
\left|\sum_{\alpha\in F(\psi_j)\atop D_1(\psi_j,\alpha)={\bf d}}c_j\right|^2\frac{1}{|d_1|},
\end{equation}
where ${\bf d}=(d_0,d_1)$.

There is a higher order version of this lower bound which works provided that for the specific $\alpha$ in $ F(\psi)$, $\psi$ 
is analytically continuable across $\partial \D$ at $\alpha$ and $|\psi(e^{i\theta})|<1$
for $e^{i\theta}$ near to, but not equal to, $\alpha$.  More detail can be found in \cite{km}, where a
somewhat larger class of maps is considered.  For such $\alpha$, the curve $\psi(
e^{i\theta})$, $e^{i\theta}$ near $\alpha$, automatically has positive and even order
of contact $2m$ with $\partial \D$ when $e^{i\theta}=\alpha$; that is,
$$\frac{1-|\psi(e^{i\theta})|^2}{|\psi(\alpha)-\psi(e^{i\theta})|^{2m}}$$
is bounded above and away from $0$ for $e^{i\theta}$ near $\alpha$.
For $k\geq 1$ the $k^{th}$ order data vector
$$D_k(\psi,\alpha)\equiv(\psi(\alpha),\psi'(\alpha),\cdots,\psi^{(k)}(\alpha))$$
makes sense. Given a finite collection $\psi_1\cdots,\psi_n$ of such maps and $k\geq 2$,
we write ${\mathcal M}_k(\alpha)$ for the set of integers $j$, $1\leq j\leq n$, for
which $F(\psi_j)$ contains $\alpha$ and the order of contact of $\psi_j$ at $\alpha$
is at least $k$.  Further, put
$${\mathcal D}_k(\alpha)=\{D_k(\psi_j,\alpha):j\in{\mathcal M}_k(\alpha)\},$$
the set of possible $k^{th}$ order data vectors at $\alpha$ for associated
orders of contact at least $k$.  With this notation, Theorem 5.7 in \cite{km}
states that for any $k\geq 2$ and complex constants $c_1\cdots, c_n$,
\begin{equation}\label{lb2}
\|c_1C_{\psi_1}+\cdots+c_nC_{\psi_n}\|_e^2\geq \sum_{{\bf d}\in{\mathcal D}_{k-1}(\alpha)}
\left|\sum_{j\in{\mathcal M}_k(\alpha)\atop D_{k-1}(\psi_j,\alpha)={\bf d}}c_j\right|^2\frac{1}{|d_1|},
\end{equation}
where ${\bf d}=(d_0,d_1,\cdots,d_{k-1})$.

For the case $k=2$, we need a calculation which 
appears in the proof of a more delicate version of the inequality~(\ref{lb2});
see Lemma 5.9 in \cite{km}.  Given $\psi$ as above, with $\alpha$ in $ F(\psi)$, convert
it into a self-map $u$ of the upper half-plane fixing the origin by $u=\tau_{\psi(\alpha)}
\circ\psi\circ \tau_{\alpha}^{-1}$.  Given our finite collection $\psi_1\cdots,\psi_n$,
associate $u_j$ to $\psi_j$ in this manner.  For $D>0$ we write $\Gamma_{\alpha,D}$
for the locus of the equation $(1-|z|^2)/(|\alpha-z|^2)=4D$, a circle internally tangent
to $\partial \D$ at $\alpha$.  We have
\begin{equation}\label{lb3}
\lim_{\Gamma_{\alpha,D}}\left\|(\overline{c_1}C^*_{\psi_1}+\cdots+\overline{c_n}C^*_{\psi_n})
\frac{k_z}{\|k_z\|}\right\|^2= \sum_{{\bf d}\in{\mathcal D}_1(\alpha)}
\left\|\sum_{j\in{\mathcal M}_2(\alpha)\atop D_1(\psi_j,\alpha)={\bf d}}
\overline{c_j}k^+_{w_j}\right\|^2
_{H^2_+}
\end{equation}
where $k_z$ is the Szego kernel for the Hardy space $H^2$ in the disk,
$H^2_+$ is the Hardy space of the right half-plane $\Omega$, $k^+_w(z)=1/(\overline{w}+z)$
is its reproducing kernel,  $w_j=u_j'(0)/2-iDu_j''(0)$, and the limit is taken
as $z\rightarrow \alpha$ along $\Gamma_{\alpha,D}$.
Since $u''_j(0)$ necessarily has non-negative imaginary part, $w_j$ is a complex number automatically lying in $\Omega$.  
For further discussion of this circle of
ideas, see \cite{km}.  We note for future reference that a non-automorphism linear-fractional
self-map $\psi$ of $\D$ has order of contact two at the unique point in $F(\psi)$.

Finally, we need a variant of a result of Berkson \cite{b} and Shapiro and Sundberg \cite{ss}, 
which states that if $\psi_1,\cdots,\psi_n$ are distinct
analytic self-maps of $\D$ and $J(\psi_i)=\{e^{i\theta}:|\psi_i(e^{i\theta})|=1\}$,
then for any complex constants $c_1,\cdots,c_n$,
\begin{equation}\label{lbext}
\|c_1C_{\psi_1}+\cdots+c_nC_{\psi_n}\|_e^2\geq \frac{1}{2\pi}\sum_{j=1}^{n}|c_j|^2|J(\psi_j)|,
\end{equation}
where $|J|$ denotes the Lebesgue measure of $J$; see Exercise 9.3.2 in \cite{cmbook}.

\begin{thm}\label{compacts}
Let $\psi$ be an analytic self-map of $\D$ such that $C_{\psi}$ lies in
$C^*(\cp,{\mathcal K})$, where $\varphi$ is as in~(\ref{defofph}).  If $F(\psi)$ is empty,
then $C_{\psi}$ is compact.
\end{thm} 

\begin{proof}
Suppose that $C_{\psi}$ lies in $C^*(\cp,{\mathcal K})$ and $F(\psi)$ is empty.
We want to show that $C_{\psi}$ is compact, or equivalently, that the matrix
function
$$\Psi(C_{\psi})=\left[\begin{array}{lr}
f_2&f_3\\f_4&f_1\end{array}\right]$$
is identically zero on $[0,s]$.
Given a small $\epsilon>0$ (size to be specified later), there exists $B$ in ${\mathcal P}$
given by Equation~(\ref{defofB}) and equivalently by Equation~(\ref{newB}), such that $\|C_{\psi}-B\|<\epsilon$.  
If we write
$$Y_1=f(\cp^*\cp), Y_2=g(\cp\cp^*), Y_3=\cp p(\cp^*\cp), Y_4=\cp^*q(\cp\cp^*),$$
and $Y=Y_1+Y_2+Y_3+Y_4$, it is clear that $A_k\equiv Y_k\ (\mbox{mod }{\mathcal K})$ for 
each $i$, and $\Psi(Y)=\Psi(A)$, where $A=A_1+A_2+A_3+A_4$.  Now using the 
representation (\ref{newB}) for $B$, we have
$$\|C_{\psi}-cC_z-A\|_e<\epsilon.$$
Since $A$ is a finite linear combination of composition operators, we see from 
the inequality~(\ref{lbext}) that
$$\epsilon^2>\|C_{\psi}-cC_z-A\|_e^2\geq\frac{|J(\psi)|}{2\pi}+|c|^2.$$
From this we find that $|c|<\epsilon$, and since $\epsilon>0$ is arbitrary,
$|J(\psi)|=0$.  In particular we have $\|C_{\psi}-A\|_e<2\epsilon$, hence
$$\left\|\left[\begin{array}{lr}
f_2-g&f_3-rp\\f_4-rq&f_1-f\end{array}\right]\right\|=\left\|\Psi(C_{\psi}-Y)\right\|=
\|C_{\psi}-A\|_e<2\epsilon.$$
It follows that $|f_3(t)-\sqrt{t}p(t)|<2\epsilon$ for $0\leq t\leq s$, and
similarly for the other three matrix entries.

We will show that $f_3$ vanishes identically on $[0,s]$.  Suppose not, so
that its supremum norm $\|f_3\|_{\infty}$ is positive.  Without loss
of generality we may assume $8\epsilon<\|f_3\|_{\infty}$.  It follows that
there is a non-degenerate closed subinterval $I$ of $[0,s]$, depending only on
$f_3$ and not containing zero,
with $\sqrt{t}|p(t)|\geq\|f_3\|_{\infty}/2$ for $t$ in $I$. Thus
\begin{equation}\label{est1}
\int_I|p(t)|^2dt \geq \frac{\|f_3\|^2_{\infty}|I|}{4s}.
\end{equation}
We return to this inequality below.

Now we want to apply Equation~(\ref{lb3}) to the linear combination $A$,
which we write as 
\begin{equation}\label{A}
A=c_1C_{\psi_1}+\cdots c_mC_{\psi_m}.
\end{equation}
Recall that the normalized Szego kernel functions $k_z/\|k_z\|$
tend
to zero weakly as $|z|\rightarrow 1$, and so
$$\|T\|_e=\|T^*\|_e\geq\limsup_{|z|\rightarrow 1}\left\|T^*\left(\frac{k_z}{\|k_z\|}
\right)\right\|$$
for any bounded operator $T$ on $H^2$.  
The linear-fractional maps $\psi_1,\cdots,\psi_m$ in Equation~(\ref{A})
are taken from the four lists in~(\ref{lists}).  The maps in each of
these lists have a common angular derivative set (a singleton) and a single
common first order data vector.  For example,
the maps $\psi_i$ from the first list all have $F(\psi_i)=\{\eta\}$ and
first order data vector $D_1(\psi_i,\eta)=(\eta,1)$, which we call ${\bf d}_1$.  
The following table summarizes the corresponding information for each of
the four lists:
\begin{center}
\renewcommand{\arraystretch}{1.45}
\framebox{
\begin{tabular}{l|c|l}
\multicolumn{3}{c}{\sc Table I}\\ \hline\hline
\multicolumn{1}{c}{} &
\multicolumn{1}{|c|}{} &
\multicolumn{1}{c}{unique first-order} \\[-3mm]
\multicolumn{1}{c}{$\psi_i$ chosen from} &
\multicolumn{1}{|c|}{$F(\psi_i$)} &
\multicolumn{1}{c}{data vector} \\ \hline
$\{(\varphi \circ \sigma)_{n_{1}}: \ n_1 \geq 1\}$ & $\{\eta\}$ & \quad ${\bf d}_1 = (\eta,1)$
\\ \hline
$\{(\sigma \circ \varphi)_{n_{2}}: \ n_2 \geq 1\}$ & $\{\zeta\}$ & \quad${\bf d}_2 = (\zeta,1)$
\\ \hline
$\{(\varphi \circ \sigma)_{n_{3}}\circ\varphi: \ n_3 \geq 0\}$ & $\{\zeta\}$ & \quad ${\bf d}_3 = (\eta,\varphi'(\zeta))$
\\ \hline
$\{(\sigma \circ\varphi)_{n_{4}}\circ\sigma: \ n_4 \geq 0\}$ & $\{\eta\}$ & \quad ${\bf d}_4 = (\zeta,\sigma'(\eta))$
\end{tabular}
}
\end{center}
Since $C_{\psi}^*(k_z)=k_{\psi(z)}$, our hypothesis that $F(\psi)$ is empty says
exactly that 
$$\lim_{|z|\rightarrow 1}\left\|C_{\psi}^*\frac{k_z}{\|k_z\|}\right\|^2=
\lim_{|z|\rightarrow 1}\frac{1-|z|^2}{1-|\psi(z)|^2}=0$$
(see \cite{cmbook}, p.132) and thus
\begin{eqnarray*}
4\epsilon^2>\|C_{\psi}-A\|_e^2&\geq &\limsup_{|z|\rightarrow 1}\left\|(C_{\psi}^*-A^*)
\frac{k_z}{\|k_z\|}\right\|^2\\
&\geq&\lim_{\Gamma_{\zeta,D}}\left\|(\overline{c_1}C_{\psi_1}^*+\cdots+\overline{c_m}
C_{\psi_m}^*)\frac{k_z}{\|k_z\|}\right\|^2.
\end{eqnarray*}
We evaluate the limit on the right via Equation~(\ref{lb3}) with $\alpha=\zeta$.
Note that ${\mathcal D}_1(\zeta)=\{{\bf d}_2,{\bf d}_3\}$.  Discarding the (necessarily
non-negative) ${\bf d}_2$ term from the right-hand side of Equation~(\ref{lb3}) yields
\begin{equation}\label{12}
4\epsilon^2>\lim_{\Gamma_{\zeta,D}}
\left\|(\overline{c_1}C_{\psi_1}^*+\cdots+\overline{c_m}C_{\psi_m}^*)
\frac{k_z}{\|k_z\|}\right\|^2\geq\left\|\sum_{\zeta\in F(\psi_i)\atop D_1(\psi_i,\zeta)={\bf d}_3}\overline{c_i}
k_{w_i}^+\right\|^2_{H^2_+}.
\end{equation}
We can relabel $\psi_1,\cdots,\psi_m$ so that the relevant $\psi_i's$ occur at the beginning, starting with
$i=0$.  Then the right-hand side of (\ref{12}) becomes
\begin{eqnarray}\label{13}
\left\|\sum_{k=0}^{n}\overline{c_k}k_{w_k}^+\right\|^2_{H^2_+}=\sum_{i,j=0}^n
\overline{c_i}c_j\frac{1}{\overline{w_i}
+w_j}
&=&\sum_{i,j=0}^n\overline{c_i}c_j\int_0^1t^{\overline{w_i}+w_j-1}dt\\
&=&\int_0^1\left|\sum_{k=0}^nc_kt^{w_k}\right|^2t^{-1}dt\nonumber
\end{eqnarray}
for appropriate $n\leq m-1$.

Let us evaluate $c_k$ in terms of the polynomial $p$ occuring in the upper right
entry of the matrix  
function $\Psi(A)=\Psi(Y)$.  If 
$$p(z)=\sum_{k=0}^nb_kz^k,$$
we have
\begin{eqnarray*}
Y_3=\cp p(\cp^*\cp)&\equiv&\cp p(sC_{\sigma}\cp)(\mbox{mod }{\mathcal K})\\
&=&\cp p(sC_{\varphi\circ\sigma})\\
&=&\sum_{k=0}^nb_ks^kC_{(\varphi\circ\sigma)_k\circ\varphi},\\
&=&A_3,
\end{eqnarray*}
so that, relabeling if necessary, $\psi_k=(\varphi\circ\sigma)_k\circ\varphi$ and $c_k=b_ks^k$
for $k=0,1,\cdots,n$. 

Next we compute $w_k$ for $k=0,1,\cdots,n$.  Let us convert $\psi_k$ into
a self-map $U_k$ of the upper half-plane fixing the origin as described
prior to Equation~(\ref{lb3}): $U_{k}=\tau_{\eta}\circ\psi_k\circ\tau_{\zeta}^{-1}$.
We can do the same for the composition factors of $\psi_k=(\varphi\circ\sigma)_k\circ\varphi$.
The map $\varphi\circ\sigma$ is a positive parabolic non-automorphism with fixed
point $\eta$.  Let $a>0$ be its translation number, so that $\varphi\circ\sigma=\rho_a$.
Thus for $k\geq 1$, $(\varphi\circ\sigma)_k=\rho_{ka}$, and its half-plane transplant
$u_k=\tau_{\eta}\circ\rho_{ka}\circ\tau_{\eta}^{-1}$ satisfies $u_k''(0)=2ika$
by Equation~(\ref{u}).  We write $v$ for the half-plane version of $\varphi:
v=\tau_{\eta}\circ\varphi\circ\tau_{\zeta}^{-1}$.  We have
$$U_k=\tau_{\eta}\circ\psi_k\circ\tau_{\zeta}^{-1}=(\tau_{\eta}\circ\rho_{ka}\circ\tau_{\eta}^{-1})
\circ(\tau_{\eta}\circ\varphi\circ\tau_{\zeta}^{-1})=u_k\circ v.$$
Now $u_k'(0)=\rho_{ka}'(\eta)=1$ and $v'(0)=|\varphi'(\zeta)|=\frac{1}{s}$, and we find 
$$U_k'(0)=\frac{1}{s}\mbox{ and }U_k''(0)=v''(0)+\frac{2ika}{s^2}.$$
From the discussion following Equation~(\ref{lb3}) we see that
$$w_k=\frac{1}{2s}-iDv''(0)+k\left(\frac{2aD}{s^2}\right).$$
To this point $D$ has been an arbitrary positive number. 
Let us choose $D$ so that $2aD/s^2=1$ and put $\mu=\frac{1}{2s}-iDv''(0)$, a complex
number with positive real part.  Thus $w_k=\mu + k$ and we can write the right hand
side of Equation~(\ref{13}) as
\begin{eqnarray*}
\int_0^1\left|\sum_{k=0}^n b_ks^kt^k\right|^2t^{2\mbox{Re}\mu-1}dt&=&\int_0^1|p(st)|^2t^{2\mbox{Re}\mu-1}dt\\
&=&\frac{1}{s^{2\mbox{Re}\mu}}\int_0^s|p(t)|^2t^{2\mbox{Re}\mu-1}dt
\end{eqnarray*}
We consider two cases:  if $2\mbox{Re}\mu-1\geq 0$ then this last integral is at least
$$\frac{t_0^{2\mbox{Re}\mu-1}}{s^{2\mbox{Re}\mu}}\int_I|p(t)|^2dt$$
where $t_0>0$ is the left-hand endpoint of $I$, and if
$2\mbox{Re}\mu-1<0$ it is at least 
$$\frac{1}{s}\int_I|p(t)|^2dt.$$  For small enough $\epsilon>0$, either case of
this inequality, combined with the inequality (\ref{12}) and Equation~(\ref{13}),
is incompatible with the inequality~(\ref{est1}), yielding the desired contradiction.
It follows that $f_3\equiv 0$ on $[0,s]$.  Entirely similar arguments show that
$f_1,f_2$ and $f_4$ vanish identically on $[0,s]$, hence $\Psi(C_{\psi})=0$.
\end{proof}

With Theorem~\ref{compacts} in hand, we can present our necessary conditions for membership
in $C^*(C_{\varphi},{\mathcal K})$.  
\begin{thm}\label{necessity}
Let $\varphi$ be as in (\ref{defofph}).
Suppose $\psi$ is an analytic self-map of $\D$ with $C_{\psi}$ lying in $C^*(\cp,{\mathcal K})$
and $C_{\psi}$ not compact.
Then either $\psi(z)=z$ or one of the following holds:
\begin{itemize}
\item[(a)] $F(\psi)=\{\zeta\}, \psi(\zeta)=\eta$ and $\psi'(\zeta)=\varphi'(\zeta)$.
\item[(b)] $F(\psi)=\{\zeta\},\psi(\zeta)=\zeta$ and $\psi'(\zeta)=1$.
\item[(c)] $F(\psi)=\{\eta\}, \psi(\eta)=\zeta$ and $\psi'(\eta)=1/\varphi'(\zeta)$.
\item[(d)] $F(\psi)=\{\eta\}, \psi(\eta)=\eta$ and $\psi'(\eta)=1$.
\item[(e)] $F(\psi)=\{\zeta,\eta\}$ with $\psi(\zeta)=\eta,\psi'(\zeta)=\varphi'(\zeta),
\psi(\eta)=\eta$ and $\psi'(\eta)=1$.
\item[(f)] $F(\psi)=\{\zeta,\eta\}$ with $\psi(\eta)=\zeta, \psi'(\eta)=1/\varphi'(\zeta), 
\psi(\zeta)=\zeta$
and
$\psi'(\zeta)=1$.
\end{itemize}
\end{thm}
\begin{proof}
If $\psi$ has no finite angular derivative, then Theorem~\ref{compacts} guarantees that $C_{\psi}$ is
compact. Thus we may assume $F(\psi)$ is non-empty.  We also assume $\psi$ is not the identity, else
there is nothing to prove.
If $C_{\psi}$ is in $C^*(C_{\varphi},{\mathcal K})$, the density of the polynomial subalgebra
${\mathcal P}$ says that given $\epsilon$, we may find a scalar $c$ and a finite linear
combination $A$ of composition operators with associated maps from  the lists (\ref{lists})
so that
$$
\|C_{\psi}-A-cC_z\|_e<\epsilon.
$$
As in the beginning of the proof of Theorem~\ref{compacts}, we may then conclude that 
$|\psi(e^{i\theta})|<1$ a.e., $|c|<\epsilon$, and that
\begin{equation}\label{approx}
\|C_{\psi}-A\|<2\epsilon.
\end{equation}
The self-maps of $\D$ which define the composition operators in the linear
combination $A$ appear among those in Table I above, along with their angular
derivative sets (all singletons) and first order data vectors.  Suppose that
$\lambda$ is in $F(\psi)$ and $D_1(\psi,\lambda)={\bf d}$.  If 
$\lambda$ is not in $\{\zeta,\eta\}$, then the inequality~(\ref{lb1}) gives
\begin{equation}\label{toobig}
\|C_{\psi}-A\|_e^2\geq\frac{1}{|\psi'(\lambda)|},
\end{equation}
contradicting (\ref{approx}).  Similarly, if $\lambda=\zeta$ and ${\bf d}$
is neither ${\bf d_2}$ nor ${\bf d_3}$ from Table I, or if 
$\lambda=\eta$ and ${\bf d}$ is neither ${\bf d_1}$ nor ${\bf d_4}$, the inequality
(\ref{lb1}) and Table I again imply (\ref{toobig}).  It follows
that if $F(\psi)$ is a singleton, one of conditions (a)-(d) must hold.

The remainder of the proof considers the possibility that $F(\psi)=\{\zeta,\eta\}$.
The Julia-Caratheodory theory says a non-identity analytic self-map
of $\D$ cannot have fixed points at distinct points $\zeta, \eta$ in $\partial \D$
with derivative $1$ at each. If we have both $\psi(\zeta)=\eta, \psi'(\zeta)=\varphi'(\zeta)$
and $\psi(\eta)=\zeta,\psi'(\eta)=1/\varphi'(\zeta)$, then
$\psi\circ\psi$ fixes both $\zeta$ and $\eta$ with derivative $1$ at each, so 
that $\psi\circ\psi$ is the identity map,
contradicting the observation above that $|\psi(e^{i\theta})|<1$ almost everywhere. 
Thus if $F(\psi)=\{\zeta,\eta\}$, either (e) or (f) must hold, completing the proof.
\end{proof}

We will see in Section 5 that there are indeed maps $\psi$ of types (e) and (f) in Theorem~\ref{necessity}
for which $C_{\psi}$ belongs to $C^*(\cp,{\mathcal K})$.

\section{The $C^*$-algebra induced by parabolic non-automorphisms}

Let us write ${\mathcal B}(H^2)$ for the algebra of bounded operators on $H^2$.  
A bounded operator $T$ on $H^2$ is essentially normal if $T$ and $T^*$ commute modulo
${\mathcal K}$; normal operators and compact operators give trivial
examples of essentially normal operators.  The only normal composition operators $C_{\psi}$ are those
of the form $\psi(z)=az, |a|\leq 1$.  Bourdon, Levi, Narayan, and Shapiro \cite{blns} showed
that if $\psi$ is linear-fractional with $\|\psi\|_{\infty}=1$ and not a rotation,
then $C_{\psi}$ is essentially normal exactly when $\psi$ is a parabolic non-automorphism.
Let us select a point $\gamma$ in $\partial \D$ and consider the set $\{\rho_a:\mbox{Re }a>0\}$ of all
parabolic non-automorphisms fixing $\gamma$.  Here, as earlier, $a$ is the translation number
for $\rho_a$.  Any two of the maps $\rho_a$ commute under composition and in fact $\rho_a\circ\rho_b=
\rho_{a+b}$, so $C_{\rho_a}$ and $C_{\rho_b}$ commute.  One can easily
check that the Krein adjoint of $\rho_a$ is $\rho_{\overline{a}}$.  Since $\rho_a'(\gamma)=1$,
it follows from Equation~(\ref{adjformula}) that $C_{\rho_a}^*=C_{\rho_{\overline{a}}}+K$ for some compact operator
$K$.  A recent theorem of Montes-Rodr\'iguez, Ponce-Escudero and Shkarin \cite{montes}
shows that $C_{\rho_a}$ is irreducible.  
Moreover $C^*(C_{\rho_a})$, the unital $C^*$-algebra generated by $C_{\rho_a}$, 
contains the commutator of $C_{\rho_a}$ and $C_{\rho_a}^*$ which we
know is compact but non-zero.  Thus $C^*(C_{\rho_a})$ contains ${\mathcal K}$
and $C^*(C_{\rho_a})/{\mathcal K}$ is commutative.
Now let ${\mathbb P}_{\gamma}$ denote the set of all composition operators $C_{\rho}$,
where $\rho$, fixing $\gamma$, ranges over $\{\rho_a:\mbox{Re }a>0\}$.  We write $C^*({\mathbb P}_{\gamma})$ for 
the unital $C^*$-algebra generated by the operators in ${\mathbb P}_{\gamma}$.  Clearly $C^*({\mathbb P}_{\gamma})$
contains ${\mathcal K}$, and, by the above remarks, $C^*({\mathbb P}_{\gamma})/{\mathcal K}$ is
also commutative.  In this section we compute and apply the Gelfand representation
of this quotient algebra.

We begin with two lemmas.

\begin{lemma}\label{positive}
For $a>0$ there is an operator $A\geq 0$ and a compact
operator $K$ with 
$C_{\rho_a}=A+K.$
\end{lemma}
\begin{proof}
Then 
\begin{eqnarray*}
C_{\rho_{a/2}}&=&\frac{1}{2}(C_{\rho_{a/2}}+C_{\rho_{a/2}}^*)+\frac{1}{2}
(C_{\rho_{a/2}}-C_{\rho_{a/2}}^*)\\
&\equiv&B+J
\end{eqnarray*}
where $B$ is self-adjoint and $J$ is compact.
Thus $C_{\rho_a}=C_{\rho_{a/2}}C_{\rho_{a/2}}=(B+J)^2=
B^2+(BJ+JB+J^2)$.  Since $B^2$ is positive and $BJ+JB+J^2$ is compact,
we are done.
\end{proof}

\begin{lemma}\label{dashcond}
Let $a,b$ be positive with $b/a=m/n$, with $m$ and $n$ positive integers.
Suppose $0<\lambda\leq 1$ and there is a sequence $f_k$ of unit vectors in $H^2$
converging weakly to zero such that
$$\|(C_{\rho_a}-\lambda)f_k\|\rightarrow 0.$$
Then
$$\|(C_{\rho_b}-\lambda^{m/n})f_k\|\rightarrow 0.$$
\end{lemma}
\begin{proof}
First observe that
$$C_{\rho_a}^m-\lambda ^m=[C_{\rho_a}^{m-1}+\lambda C_{\rho_a}^{m-2}+\cdots
+\lambda^{m-2}C_{\rho_a}+\lambda^{m-1}][C_{\rho_a}-\lambda].$$
In particular, $\|(C_{\rho_a}^m-\lambda^m)f_k\|\rightarrow 0$ as 
$k\rightarrow\infty$.  Since $C_{\rho_a}^m=C_{\rho_b}^n,$
$$\|(C_{\rho_b}^n-\lambda^m)f_k\|\rightarrow 0.$$
Also note that we may factor $
C_{\rho_b}^n-\lambda ^m=C_{\rho_b}^n-(\lambda^{m/n})^n$ as
$$
[C_{\rho_b}^{n-1}+\lambda^{m/n} C_{\rho_b}^{n-2}+\cdots
+(\lambda^{m/n})^{n-2}C_{\rho_b}+(\lambda^{m/n})^{n-1}]
[C_{\rho_b}-\lambda^{m/n}].$$
Apply Lemma~\ref{positive} to $C_{\rho_b}$ to write
$$C_{\rho_b}^{n-1}+\lambda^{m/n} C_{\rho_b}^{n-2}+\cdots
+(\lambda^{m/n})^{n-2}C_{\rho_b}+(\lambda^{m/n})^{n-1}I=T+(\lambda^{m/n})^{n-1}I+K$$
for some positive $T$ and compact $K$.
We have 
\begin{equation}\label{1}
\|(C_{\rho_b}^n-\lambda^m)f_k\|=\|(T+(\lambda^{m/n})^{n-1}+K)(C_{\rho_b}-
\lambda^{m/n})f_k\|.
\end{equation}
Since $K$ is compact, $\|K(C_{\rho_b}-\lambda^{m/n})f_k\|\rightarrow 0$
as $k\rightarrow \infty$.  Since the left-hand side of Equation~(\ref{1})
goes to $0$ as $k\rightarrow \infty$, we see, writing $c=(\lambda^{m/n})^{n-1}$
that
$$\|(T+cI)(C_{\rho_b}-\lambda^{m/n})f_k\|\rightarrow 0.$$
But
$\|(T+cI)(C_{\rho_b}-\lambda^{m/n})f_k\|^2$ is equal to
\begin{eqnarray*}
\|T(C_{\rho_b}-\lambda^{m/n})f_k\|^2
&+&2c\langle T(C_{\rho_b}-\lambda^{m/n})f_k,(C_{\rho_b}-\lambda^{m/n})f_k\rangle\\
&+&c^2\|(C_{\rho_b}-\lambda^{m/n})f_k\|^2\geq
c^2\|(C_{\rho_b}-\lambda^{m/n})f_k\|^2
\end{eqnarray*}
where the last inequality uses the positivity of $T$,
so that $\|(C_{\rho_b}-\lambda^{m/n})f_k\|\rightarrow 0$.
\end{proof}

The essential spectrum $\sigma_e(T)$ of a bounded operator on $H^2$
is by definition the spectrum of the coset $[T]$ in ${\mathcal B}(H^2)/{\mathcal K}$.
We recall from \cite{kmm} that if $a>0, \sigma_e(C_{\rho_a})=[0,1]$.
We will need the notion of joint essential spectrum, which is treated by Dash in \cite{dash}.  
If $\mbox{Re }a>0$, the coset $[C_{\rho_a}]$ of $C_{\rho_a}$ modulo ${\mathcal K}$ will also be denoted
by $x_a$.  By either Lemma~\ref{positive} or the discussion preceeding it,
$x_a$ is self-adjoint.  Given $a$ and $b$, the joint essential spectrum $\sigma_e(C_{\rho_a},C_{\rho_b})$ is defined
to be the joint spectrum $\sigma(x_a,x_b)$ of the pair $x_a,x_b$ in the Calkin algebra
${\mathcal B}(H^2)/{\mathcal K}$.  This set coincides with the joint spectrum in the commutative
unital subalgebra $C^*(x_a,x_b)$ generated by $x_a$ and $x_b$.  If ${\mathcal M}$ is the maximal ideal
space of this algebra, and $\widehat{\rule{4pt}{0pt}}$ denotes the Gelfand transform, then the map 
$\ell\mapsto (\widehat{x_a}(\ell),\widehat{x_b}(\ell))$ is a 
homeomorphism of ${\mathcal M}$ onto $\sigma(x_a,x_b)$.  Let us assume that $a$ and $b$ are positive.
A theorem of Dash \cite{dash} states, in this context, 
using $C_{\rho_a}^*\equiv C_{\rho_a}(\mbox{mod }{\mathcal K})$ and similarly for $C_{\rho_b}$, that
$(\lambda,\mu)$ lies in $\sigma_e(C_{\rho_a},C_{\rho_b})$ if and only if there exists a sequence
$\{f_k\}$ of unit vectors in $H^2$, converging weakly to zero, 
such that $\|(C_{\rho_a}-\lambda)f_k\|$ and $\|(C_{\rho_b}-\mu)f_k\|$
both tend to zero as $k\rightarrow \infty$.
\begin{cor}\label{jointspec}
Suppose that $a,b$ are positive and $b/a$ is rational.  Then
$$\sigma_e(C_{\rho_a},C_{\rho_b})=\{(t^a,t^b):0\leq t\leq 1\}.$$
\end{cor}
\begin{proof}
We know that $$\sigma_e(C_{\rho_a},C_{\rho_b})\subset
\sigma_e(C_{\rho_a})\times\sigma_e(C_{\rho_b})=
[0,1]\times[0,1].$$
Let $0<\lambda\leq 1$.  Since $\lambda\in\sigma_e(C_{\rho_a})$,
we may find unit vectors $f_k$ with $f_k\rightarrow 0$ weakly and
$\|(C_{\rho_a}-\lambda)f_k\|\rightarrow 0$. By Lemma~\ref{dashcond}, 
$\|(C_{\rho_b}-\lambda^{b/a})f_k\|\rightarrow 0$
as well, so that $(\lambda,\lambda^{b/a})$ is in $\sigma_e
(C_{\rho_a},C_{\rho_b})$ by Dash's theorem.  Setting
$t^a=\lambda$ we have $\{(t^a,t^b):0<t\leq 1\}\subset
\sigma_e(C_{\rho_a}, C_{\rho_b})$.  The set on the 
right is compact in $\mathbb{R}^2$, so it contains $(0,0)$ as well,
giving $\{(t^a,t^b):0\leq t\leq 1\}\subset\sigma_e
(C_{\rho_a},C_{\rho_b})$.  Conversely, if $(\lambda,\mu)
\in\sigma_e(C_{\rho_a},C_{\rho_b})$, Dash's theorem
gives the existence of a sequence of unit vectors $f_k$ converging
weakly to $0$ with both
$\|(C_{\rho_a}-\lambda)f_k\|$ and
$\|(C_{\rho_b}-\mu)f_k\|$ tending to zero as $k\rightarrow\infty$.
If $\lambda>0$, $\|(C_{\rho_b}-\lambda^{b/a})f_k\|\rightarrow 0$
by Lemma~\ref{dashcond}.  Thus $\mu=\lambda^{b/a}$ and
$(\lambda,\mu)=(\lambda,\lambda^{b/a})=(t^a,t^b)$
for some $t$, $0\leq t\leq 1$.  If $\mu>0$, the symmetric
result says $\lambda=\mu^{a/b}>0$, and, putting
$\mu=t^b$, $(\lambda,\mu)=(\mu^{a/b},\mu)=(t^a,t^b)$,
again of the desired form.  As for $(0,0)$, we already
know it lies in $\sigma_e(C_{\rho_a},C_{\rho_b})$,
and of course it has the form $(0^a,0^b)$.
\end{proof}

We will need the fact that on the domain $\Omega=\{a:\mbox{Re }a>0\}$
the map $a\mapsto C_{\rho_a}$ is a holomorphic function of $a$ in the operator norm
topology; see for example the discussion in the proof of 
Theorem 6.1 in \cite{cowen}.  We continue to denote 
the coset of $C_{\rho_a}$ by $x_a$ and to keep in mind that
when $a>0$, $\sigma(x_a)=[0,1]$.

\begin{thm}\label{Gamma}
There is a unique $*$-isomorphism $\Gamma:C([0,1])\rightarrow
C^*({\mathbb P}_{\gamma})/{\mathcal K}$ such that $\Gamma(t^a)=[C_{\rho_a}]$
for all $a\in\Omega$.
\end{thm}
\begin{proof}
First consider $a=1$ and $x_1=[C_{\rho_1}]$. 
Since $\sigma(x_1)=[0,1]$ we may define a $*$-isomorphism $\Gamma:C([0,1])\mapsto
C^*(x_1)$ by sending $p$ to $p(x_1)$ for any polynomial $p$.  Fix any rational number
$r>0$.
By Corollary ~\ref{jointspec},
$$\sigma(x_1,x_r)=\sigma_e(C_{\rho_1},C_{\rho_r})=\{(t,t^r):0\leq t\leq 1\}.$$
The map $p(x_1,x_r)\mapsto p(z_1,z_2)$, where $p$ is a two-variable polynomial, 
extends to a unique $*$-isomorphism of
$C^*(x_1,x_r)$ onto $C(\sigma(x_1,x_r))$.  Since $\sigma(x_1,x_r)$ is homeomorphic to $[0,1]$ via the map
$t\mapsto (t,t^r)$, we see that $p(x_1,x_r)\mapsto p(t,t^r)$ defines a 
$*$-isomorphism of $C^*(x_1,x_r)$ onto $C([0,1])$.  Let $\tilde{\Gamma}$ denote
the inverse of this map, that is $\tilde{\Gamma}:p(t,t^r)\mapsto p(x_1,x_r)$.  Since
polynomials in $t$ span $C([0,1])$, $x_1$ generates the $C^*$-algebra $C^*(x_1,x_r)$ and
$C^*(x_1)=C^*(x_1,x_r)$.  It follows that $\tilde{\Gamma}=\Gamma$. Since $r$ is 
arbitrary in the set ${\mathbb Q}_+$ of positive rationals, we have shown that
$$C^*(\{x_r:r\in {\mathbb Q}_+\})=C^*(x_1).$$  Moreover,
$\Gamma(t^r)=x_r, r\in {\mathbb Q}_+$.
It is easy to see that the map $a\mapsto t^a$ is a norm-holomorphic map
of the right half plane into $C([0,1])$ and thus that $a\mapsto \Gamma(t^a)$ is norm-
holomorphic from the right half-lane to ${\mathcal B}(H^2)/{\mathcal K}$,
as is the function $a\mapsto x_a$.  We have seen that these functions agree
on ${\mathbb Q}_+$, hence they must agree on the right half-plane $\Omega$.
\end{proof}

We record three immediate consequences.

\begin{cor}
If $a_1,\cdots,a_n$ lie in the right half-plane $\Omega$, then
$$\sigma_e(C_{\rho_{a_1}}, \cdots,C_{\rho_{a_n}})=\{(t^{a_1},\cdots,t^{a_n}):0\leq t\leq 1\}.$$
\end{cor}

\begin{cor}
If $\rho$ is a parabolic non-automorphism fixing $\gamma$, then 
$C^*(C_{\rho})=C^*({\mathbb P}_{\gamma})$.
\end{cor}

\begin{cor} \label{parabolicsin}
If $\varphi$ is as in (\ref{defofph}), then
${\mathbb P}_{\zeta}$ and ${\mathbb P}_{\eta}$ are
both subsets of $C^*(\cp,{\mathcal K})$.
\end{cor}

\section{Linear-fractional maps}

The goal of this section is to find all linear-fractional $\psi$ with
$C_{\psi}$ in $C^*(\cp, {\mathcal K})$, where $\varphi$ satisfies
the conditions of ~(\ref{defofph}).  Since 
$C_{\psi}$ is compact if $\|\psi\|_{\infty}<1$, our interest is in the case
$\|\psi\|_{\infty}=1$.

\begin{lemma}\label{easycontain}
If $\varphi$ is as in~(\ref{defofph}), $C^*(\cp,{\mathcal K})$ contains $C_{\psi}$
for all linear-fractional $\psi:\D\rightarrow \D$ with $\psi(\zeta)=\eta$,
$\psi'(\zeta)=\varphi'(\zeta)$ and $\psi(\D)$ properly contained in $\varphi(\D)$.
\end{lemma}
\begin{proof}
Set $\tau=\varphi^{-1}\circ\psi$, noting that the hypothesis
$\psi(D)\subset\varphi(\D)$ means that $\tau$ is well-defined.
Since this containment is proper, and $\tau'(\zeta)=1$, $\tau$ is
a parabolic non-automorphism with fixed point $\zeta$.  Since
$\varphi\circ\tau = \psi$, $C_{\psi}=C_{\tau}\cp$.  By Corollary~\ref{parabolicsin},
$C_{\tau}\in C^*(\cp,{\mathcal K})$, from which the conclusion follows.
\end{proof}

Now consider the parabolic non-automorphism $\rho=\varphi\circ\sigma$.  The unique fixed point for
$\rho$ and its iterates $(\rho)_n$ is $\eta$.  Fix an integer $n\geq 1$ and let $\varphi_1=(\rho)_n\circ\varphi$.
Clearly $\varphi_1(\D)$ is properly contained in $\varphi(\D)$.  Note that $\varphi_1(\zeta)=\eta$
and, since $\rho'(\eta)=1$, $\varphi_1'(\zeta)=\varphi'(\zeta)$.  It follows from Lemma~\ref{easycontain}
that $C^*(C_{\varphi_1},{\mathcal K})$ is contained in $C^*(\cp,{\mathcal K})$.  Let
$i$ denote the corresponding inclusion map.  Since $i({\mathcal K})={\mathcal K}$, $i$
induces a $*$-homomorphism
$$\hat{i}:C^*(C_{\varphi_1},{\mathcal K})/{\mathcal K}\rightarrow C^*(\cp,{\mathcal K})/{\mathcal K}$$
given by $\hat{i}([T])=[T]$, where $[T]$ denotes the coset, modulo ${\mathcal K}$ of the
operator $T$.  Note that $\hat{i}$ is itself an inclusion.  Also observe that the map 
$\Psi:C^*(\cp,{\mathcal K})\rightarrow {\mathcal D}$ induces a $*$-isomorphism
$\Phi:C^*(\cp,{\mathcal K})/{\mathcal K}\rightarrow {\mathcal D}$ given by
$\Phi([T])=\Psi(T)$.  Let $\Phi_1$ denote the corresponding $*$-isomorphism 
$\Phi_1:C^*(C_{\varphi_1},{\mathcal K})/{\mathcal K}
\rightarrow {\mathcal D}$. Keep in mind that $\Phi_1$ should be defined by $\Phi_1([T])=\Psi_1(T)$,
where $\Psi_1:C^*(C_{\varphi_1},{\mathcal K})\rightarrow {\mathcal D}$ is associated
to $\varphi_1$ as $\Psi$ is associated to $\varphi$.  Thus if $B$ in $C^*(C_{\varphi_1},{\mathcal K})$ is 
given by ~(\ref{defofB}), {\it but with} $\varphi$ {\it replaced by} $\varphi_1$, then
$\Psi_1(B)$ is given by ~(\ref{defofpsi}). We have a commutative diagram

\vspace{.2in}
$$
\begin{array}{ccc}
C^*(C_{\varphi_{1}},\mathcal{K})/\mathcal{K}
&
\begin{picture}(30,10)
\put(0,4){\vector(1,0){30}}
\put(12,8){$\hat{i}$}
\end{picture}
&
C^*(C_{\varphi},\mathcal{K})/\mathcal{K}
\\[2mm]
\begin{picture}(10,30)
\put(5,30){\vector(0,-1){30}}
\put(-10,12){$\Phi_1$}
\end{picture}
&&
\begin{picture}(10,30)
\put(5,30){\vector(0,-1){30}}
\put(9,12){$\Phi$}
\end{picture}
\\[2mm]
\mathcal{D}
&
\begin{picture}(30,10)
\put(-30,4){\vector(1,0){88}}
\put(9,8){$\Lambda$}
\end{picture}
&
\mathcal{D}
\end{array}
$$
where $\Lambda=\Phi\circ\hat{i}\circ\Phi_1^{-1}$.  We seek to identify $\Lambda$ explicitly.
\begin{lemma}\label{Lambda}
For any element $F$ in ${\mathcal D}$, 
\begin{equation}\label{formoflambda}
(\Lambda F)(t)=F(t^{2n+1}/s^{2n}), \ 0\leq t\leq s.
\end{equation}
\end{lemma}
\begin{proof} 
For the purposes of the proof, we use $\Lambda$ to denote the map given by formula~(\ref{formoflambda}),
and then show, with this redefinition, that it coincides with $\Phi\circ\hat{i}\circ\Phi_1^{-1}$,
that is, that $\Lambda\circ\Phi_1=\Phi\circ\hat{i}$.
Recall that $C_{\varphi}^*\equiv sC_{\sigma}\mbox{ (mod }{\mathcal K})$ so that
\begin{eqnarray*}
C_{\varphi_1}=C_{\varphi}C_{(\rho)_n}&=&C_{\varphi}(C_{\sigma}C_{\varphi})^n\\
&\equiv&\frac{1}{s^n}
C_{\varphi}(C_{\varphi}^*C_{\varphi})^n \ (\mbox{mod }{\mathcal K}),
\end{eqnarray*}  and, taking adjoints, 
$C_{\varphi_1}^*\equiv\frac{1}{s^n}C_{\varphi}^*(C_{\varphi}C_{\varphi}^*)^n$ modulo the compacts.
Calculations using these two facts show that if we write $y=[C_{\varphi_1}]$ 
and $x=[\cp]$, we have, for each non-negative integer $m$,
$$(y^*y)^m=\frac{1}{s^{2nm}}(x^*x)^{(2n+1)m},$$
$$(yy^*)^m=\frac{1}{s^{2nm}}(xx^*)^{(2n+1)m}.$$
$$y(y^*y)^m=\frac{1}{s^{(2m+1)n}}x(x^*x)^{(2n+1)m+n},$$
and
$$y^*(yy^*)^m=\frac{1}{s^{(2m+1)n}}x^*(xx^*)^{(2n+1)m+n}.$$
The left-hand sides in these four equations are elements in $C^*(C_{\varphi_1},{\mathcal K})/{\mathcal K}$,
while the right-hand sides represent the same objects as elements of $C^*(\cp,{\mathcal K})/{\mathcal K}$.  We 
first act on $y(y^*y)^m$ by $\hat{i}$, followed by $\Phi$.  We then act on
$y(y^*y)^m$ by $\Phi_1$, followed by $\Lambda$ (as defined by Equation~(\ref{formoflambda})).  As the reader can see from
the following picture, we end up with a common result, the matrix function in the lower right-hand corner.

\vspace{.2in}

$$
\begin{array}{ccc}
y(y^*y)^m
&
\begin{picture}(30,10)
\put(0,4){\vector(1,0){30}}
\put(12,8){$\hat{i}$}
\end{picture}
&
\frac{1}{s^{(2m+1)n}}x(x^*x)^{(2n+1)m+n}
\\[2mm]
\begin{picture}(10,30)
\put(5,30){\vector(0,-1){30}}
\put(-10,12){$\Phi_1$}
\end{picture}
&&
\begin{picture}(10,30)
\put(5,30){\vector(0,-1){30}}
\put(9,12){$\Phi$}
\end{picture}
\\[2mm]
\left[\begin{array}{lc}
0&\sqrt{t}\ t^m\\0&0\end{array}\right]
&
\begin{picture}(30,10)
\put(0,4){\vector(1,0){30}}
\put(9,8){$\Lambda$}
\end{picture}
&
\left[\begin{array}{lc}
0&\frac{\sqrt{t}\ t^{(2n+1)m+n}}{s^{(2m+1)n}}\\0&0\end{array}\right]
\end{array}
$$

One can check that we also arrive at common values when $\Lambda\circ\Phi_1$ and $\Phi\circ\hat{i}$ 
act on $(y^*y)^m$, and similarly
for $(yy^*)^m$ and $y^*(yy^*)^m$.  Since elements of the form $(y^*y)^m$, $(yy^*)^m$, $y(y^*y)^m$ and $y^*(yy^*)^m$, 
together with the identity, span $C^*(C_{\varphi_1},{\mathcal K})/{\mathcal K}$, we
have $\Lambda\circ\Phi_1=\Phi\circ\hat{i}$ as desired
\end{proof}

It is clear from ~(\ref{formoflambda}) that $\Lambda$ is an automorphism of ${\mathcal D}$.  It follows
that $\hat{i}$ is an isomorphism and thus that $i$ has range equal to all of $C^*(\cp,{\mathcal K})$, that
is, 
\begin{equation}\label{eqalg}
C^*(C_{\varphi_1},{\mathcal K})=C^*(\cp,{\mathcal K}).
\end{equation}
More generally, we have the following result.

\begin{thm}\label{samealgebras}
Let $\psi$ be a linear-fractional map of $\D$, not an automorphism, with
$\psi(\zeta)=\varphi(\zeta)$ and $\psi'(\zeta)=\varphi'(\zeta)$, where $\varphi$ is as in ~(\ref{defofph}).
Then $C^*(C_{\psi},{\mathcal K})=C^*(\cp, {\mathcal K})$.
\end{thm}
\begin{proof}
The circles $\varphi(\partial\D)$ and $\psi(\partial\D)$ are both internally tangent to 
$\partial \D$ at $\eta$. If $\varphi(\D)$ is a proper subset of $\psi(\D)$, then
\begin{equation}\label{containment}
C^*(\cp,{\mathcal K})\subset C^*(C_{\psi},{\mathcal K})
\end{equation}
by Lemma~\ref{easycontain}.
Suppose on the other hand that $\psi(\D)\subset\varphi(\D)$.  If $a$ is the (necessarily positive)
translation number of the parabolic map $\varphi\circ\sigma$ (so that $\varphi\circ\sigma=\rho_a$ in the
terminology of Section 2), then $(\varphi\circ\sigma)_n=\rho_{na}$,
and the radius of the disk $\rho_{na}(\D)$ shrinks to zero as $n\rightarrow\infty$.  Thus there exists
$n$ with $\rho_{na}(\D)$ properly contained in $\psi(D)$.  If $\varphi_1=(\varphi\circ\sigma)_n\circ\varphi=\rho_{na}\circ
\varphi$, then $\varphi_1(\D)$ is also properly contained in $\psi(\D)$.  Since $\varphi_1(
\zeta)=\eta=\psi(\zeta)$ and
$\varphi_1'(\zeta)=\rho_{na}'(\eta)\varphi'(\zeta)=\varphi'(\zeta)=\psi'(\zeta)$,
Lemma~\ref{easycontain} implies that $C^*(C_{\psi},{\mathcal K})$ contains
$C^*(C_{\varphi_1},{\mathcal K})$, which by ~(\ref{eqalg}) coincides with $C^*(\cp,{\mathcal K})$.
The result is that (\ref{containment}) holds, whatever the relationship between the disks
$\varphi(\D)$ and $\psi(\D)$. The statement of the theorem is symmetric in $\varphi$ and $\psi$,
so symmetry implies that the containment reverse to that in (\ref{containment}) also holds, completing the proof.
\end{proof}

\begin{thm}\label{linfrac}
Suppose $\varphi$ is as in ~(\ref{defofph}).  Let $\psi$, not the identity, be any
linear-fractional self-map of $\D$ with $\|\psi\|_{\infty}=1$. Then
$C_{\psi}$ is in $C^*(\cp,{\mathcal K})$ if and only if $\psi$ is not an automorphism and
one of the following conditions holds:
\begin{itemize}
\item[(a)] $\psi(\zeta)=\eta$ and $\psi'(\zeta)=\varphi'(\zeta)$.
\item[(b)] $\psi(\zeta)=\zeta$ and $\psi'(\zeta)=1$.
\item[(c)] $\psi(\eta)=\zeta$ and $\psi'(\eta)=1/\varphi'(\zeta)$.
\item[(d)] $\psi(\eta)=\eta$ and $\psi'(\eta)=1$.
\end{itemize}
\end{thm}

\begin{proof}
The ``only if" statement follows immediately from Theorem~\ref{necessity} and the hypothesis
that $\psi$ is linear-fractional.
Conversely, let $\psi$ be a linear-fractional map which is not an automorphism.
If $\psi$ is parabolic with fixed point at either
$\zeta$ or $\eta$, the result follows from Corollary~\ref{parabolicsin}; this handles the cases (b) and (d). 
If $\psi$ is as in (a), then we have $C_{\psi}\in C^*(\cp,{\mathcal K})$ by Theorem~\ref{samealgebras}.
Finally, if $\psi$  satisfies condition (c), then
its Krein adjoint $\sigma_{\psi}$ is a linear-fractional self-map of $\D$,
not an automorphism, which satisfies condition (a), so
that $C_{\sigma_{\psi}}\in C^*(\cp,{\mathcal K})$.
Since $C_{\psi}^*\equiv s C_{\sigma_{\psi}}$ modulo the compacts, this completes the argument.
\end{proof}

The maps satisfying (a)-(d) in Theorem~\ref{linfrac} can be described more explicitly.
Given a point $\gamma$ on $\partial \D$, let us write $\rho_{\gamma,a}$ for the unique
parabolic map fixing $\gamma$ with translation number $a$.  This will be a self-map
of $\D$ when $\mbox{Re }a\geq 0$, but when $\mbox{Re }a<0$, $\rho_{\gamma,a}$ takes $\D$
onto a larger disk, whose boundary is externally tangent to $\partial \D$ at $\gamma$. Clearly,
the linear-fractional non-automorphisms of $\D$ satisfying (b) or (d) are, respectively,
$\rho_{\zeta,a}$ or $\rho_{\eta,a}$ with $\mbox{Re }a>0$.  Imaginary $a$ gives an automorphism
of $\D$, but in this case $C^*(\cp,{\mathcal K})$ does not contain the corresponding composition
operator.  

Consider now a linear-fractional non-automorphism $\psi$ of $\D$ satisfying (a).  If $\psi(\D)
\subset \varphi(\D)$, we can define $\rho=\psi\circ\varphi^{-1}$ which fixes $\eta$ and carries $\D$ to $\D$.
Moreover, $\rho'(\eta)=\psi'(\zeta)/\varphi'(\zeta)=1$, so $\rho$ is parabolic; say $\rho=\rho_{\eta,a}$
where $\mbox{Re }a\geq 0$, and we find $\psi=\rho_{\eta,a}\circ\varphi$.  On the other hand, if
$\psi$ satisfies (a) and $\varphi(\D)$ is a proper subset of $\psi(\D)$, we put $\rho=\varphi\circ\psi^{-1}$
which is again a parabolic self-map of $\D$; this time a non-automorphism.  Thus $\rho=\rho_{\eta,a}$
with $\mbox{Re }a>0$ and we find $\psi=\rho_{\eta,a}^{-1}\circ\varphi=\rho_{\eta,-a}\circ\varphi$.
If $b$ is the unique positive number with $\rho_{\eta,b}(\D)=\varphi(\D)$, then
$\rho_{\eta,-a}\circ\varphi$ is a non-automorphism self-map of $\D$ exactly when $\mbox{Re }a<b$.
Rephrasing and summarizing, we conclude that the non-automorphisms $\psi$ of $\D$ satisfying (a)
are precisely the maps of the form $\rho_{\eta,a}\circ\varphi$, $\mbox{Re }a>-b$.  Similarly, if 
$c$ is the unique positive number with $\rho_{\zeta,c}(\D)=\sigma(\D)$, then the non-automorphisms
$\psi$ satisfying (c) are exactly the maps of the form $\psi=\rho_{\zeta,a}\circ\sigma$ with
$\mbox{Re }a>-c$.  
The next result shows that the positive translation numbers $b$ and $c$ are nicely related to
each
other, and to the translation numbers of the positive parabolic non-automorphisms
$\varphi\circ\sigma$ and $\sigma\circ\varphi$.

\begin{thm}\label{transnumb}
Let $b$ and $c$ be the unique positive numbers with
$\rho_{\eta,b}(\D)=\varphi(\D)$ and $\rho_{\zeta,c}(\D)=\sigma(\D)$,
respectively.  We have $c=|\varphi'(\zeta)|b$, and moreover,
$\varphi\circ\sigma=\rho_{\eta,2b}$ and $\sigma\circ\varphi=\rho_{\zeta,2c}$.
\end{thm}
\begin{proof}
Clearly there is no loss of generality in assuming that $\zeta=1$.  The non-affine linear-fractional
self-maps of $\D$ which send $1$ to $\eta\in\partial \D$ can be written in the form
$$\varphi(z)=\eta\ \frac{(1+s+sd)z+(d-s-sd)}{z+d}$$
where $s=|\varphi'(1)|$ and $\mbox{Re }\frac{d-1}{d+1}\geq s$ (see \cite{br}).  A computation shows
that $\varphi'(1)=\eta s$ and $\varphi''(1)=-2\eta s/(1+d)$.  The image of the unit circle
under $\varphi$ is a circle with curvature $\kappa_1=|\varphi'(1)|^{-1}\mbox{Re }[1+\varphi''(1)/\varphi'(1)]=
\frac{1}{s}\mbox{Re }[1-\frac{2}{1+d}]$.  Since 
$$\sigma(z)=\frac{\overline{\eta}(1+s+s\overline{d})z-1}{-\overline{\eta}(\overline{d}-s-s\overline{d})z+\overline{d}},$$
we find that $\sigma'(\eta)=\overline{\eta}/s$ and 
$$\sigma''(\eta)=\frac{2\overline{\eta}^2(\overline{d}-s-s\overline{d})}{s^2(1+\overline{d})},$$
so that the image of the unit circle under $\sigma$ is a circle with curvature
$$\kappa_2=s\mbox{Re }\left\{1+\eta\ \frac{2\overline{\eta}^2(\overline{d}-s-s\overline{d})}
{s^2(1+\overline{d})}\cdot\frac{s}{\overline{\eta}}\right\}=2-2\mbox{Re }\frac{1}{1+d}-s.$$
A positive parabolic non-automorphism fixing $1$ and corresponding to translation
by $a$ has the form $((2+d)z-1)/(z+d)$ where $a=-2/(d+1)$; by the above calculations
the image of the unit circle under this map has curvature $1+\mbox{Re }a$.  
Thus, if $\varphi(1)=\eta$ and $|\varphi'(1)|=s$, the unique positive value $b$
such that the curvature of $\varphi(\partial\D)$ is equal to the curvature of the image
of $\partial \D$ under the positive parabolic map which corresponds to translation by $b>0$
satisfies
$$1+b=\frac{1}{s}\ \mbox{Re }\left(1-\frac{2}{1+d}\right);$$
that is, 
$$b=\frac{1}{s}\ \mbox{Re}\left(1-\frac{2}{1+d}\right)-1.$$  Similarly, the curvature of the circle
$\sigma(\partial\D)$ is equal to the curvature of the circle which is the image of the unit
circle under the positive parabolic non-automorphism corresponding to translation by $c$
precisely when $c=1-s-2\mbox{Re}\frac{1}{1+d}$.  Thus $c=sb$.  This conclusion also holds when
$\varphi$ is an affine map, $\varphi(z)=\eta(sz+1-s)$, where the computations are easier.

For the final statement, let $\psi=\rho_{\eta,-b}\circ\varphi$, so that $\psi$ is an
automorphism of $\D$ and $\varphi=\rho_{\eta,b}\circ\psi$.  Since the Krein adjoint of
an automorphism is its inverse,
we have
$$\sigma=\sigma_{\varphi}=\sigma_{\psi}\circ\sigma_{\rho_{\eta,b}}=\psi^{-1}\circ\rho_{\eta,\overline{b}}=\psi^{-1}\circ\rho_{\eta,b}$$
and thus $\varphi\circ\sigma=\rho_{\eta,b}\circ\psi\circ\psi^{-1}\circ\rho_{\eta,b}=\rho_{\eta,2b}$.
Similarly, $\sigma\circ\varphi=\rho_{1,2c}=\rho_{\zeta,2c}$.
\end{proof}

The remarks preceeding Theorem~\ref{transnumb} express the linear-fractional maps $\psi$ with
$C_{\psi}$ belonging to $C^*(\cp,{\mathcal K})$ in terms of $\varphi,\sigma, \rho_{\eta,a}$ and
$\rho_{\zeta,a}$ for appropriate ranges of the translation numbers $a$.  We describe below a 
corresponding operator-theoretic description of $C_{\psi}$ modulo ${\mathcal K}$, in terms
of the polar factors of $\cp$ and $\cp^*$.  In \cite{kmm} it was shown that every
operator $B$ in $C^*(\cp,{\mathcal K})$ has a representation generalizing Equation~(\ref{defofB}) and
having the form $B=T+K$ with $K$ compact and
\begin{equation}\label{defofT}
T=cI+f(\cp^*\cp)+g(\cp\cp^*)+Uh(\cp^*\cp)+U^*k(\cp\cp^*),
\end{equation}
where $f$ and $h$ are continuous on $\sigma(\cp^*\cp)$, $g$ and $k$ are continuous
on $\sigma(\cp\cp^*)$, all four functions vanish at zero, and $U$ is the partial isometry
polar factor (which in this case is unitary) of $\cp$.  The restrictions of $f,g,h$ and $k$ to the interval
$[0,s]$, which coincides with both of the essential spectra $\sigma_e(\cp^*\cp)$ and $\sigma_e(\cp\cp^*)$,
are uniquely determined by $B$.  We call $T$ a {\it distinguished representative} of the coset $[B]$, and
recall from \cite{kmm} that 
$$\Psi(B)=\Psi(T)=\left[\begin{array}{lr}
c+g&h\\k&c+f\end{array}\right].$$

We start with the operator $(\cp^*\cp)^a$, defined by the self-adjoint functional calculus, where
$\mbox{Re }a>0$.  Note that 
$$\cp^*\cp\equiv s C_{\sigma}C_{\varphi}\ (\mbox{mod }{\mathcal K})=sC_{\varphi\circ\sigma}=sC_{\rho_{\eta,2b}}$$
where the last equality follows from Theorem~\ref{transnumb}.  
In Theorem~\ref{Gamma}, take $\gamma=\eta$ and consider the $*-$isomorphism $\Gamma$, here called
$\Gamma_{\eta}$ to emphasize the fixed point $\eta$.
We have
$$
[(\cp^*\cp)^a]=[\cp^*\cp]^a=s^a[C_{\rho_{\eta,2b}}]^a=s^a\Gamma_{\eta}(t^{2b})^a
=s^a\Gamma_{\eta}(t^{2ba})=s^a[C_{\rho_{\eta,2ba}}].
$$
The first and fourth equalities follow, respectively, from the facts that the coset map $B\mapsto [B]$
and $\Gamma_{\eta}$ are each $*$-homomorphisms.  Relabeling $2ba$ as $a$, we see that
$s^{-\frac{a}{2b}}(\cp^*\cp)^{\frac{a}{2b}}$ is a distinguished representative of $[C_{\rho_{\eta,a}}]$
for $\mbox{Re }a>0$.  A similar argument shows that the coset $[C_{\rho_{\zeta,a}}]$ has distinguished
representative $s^{-\frac{a}{2c}}(\cp\cp^*)^{\frac{a}{2c}}$ for $\mbox{Re }a>0$.  

Now consider $\rho_{\eta,a}\circ\varphi$, which we know to be a self-map of $\D$, but not
an automorphism, when $\mbox{Re }a>-b$.  First we look at the case $\mbox{Re }a>0$.  We have
$$
C_{\rho_{\eta,a}\circ\varphi}=\cp C_{\rho_{\eta,a}}=U(\cp^*\cp)^{\frac{1}{2}}C_{\rho_{\eta,a}}
\equiv s^{-\frac{a}{2b}}U(\cp^*\cp)^{\frac{1}{2}+\frac{a}{2b}}\ (\mbox{mod }{\mathcal K})
$$
by our above discussion.  By the spectral theorem, $(C_{\varphi}^*C_{\varphi})^z$
is holomorphic for $\mbox{Re }z>0$ in the weak operator topology,
and therefore in the operator norm topology; see \cite{HP}, Theorem 3.10.1.
Thus the cosets $[C_{\rho_{\eta,a}\circ\varphi}]$ and $[s^{-\frac{a}{2b}}U(\cp^*\cp)^
{\frac{1}{2}+\frac{a}{2b}}]$ are both holomorphic ${\mathcal B}(H^2)/{\mathcal K}$-valued functions of
$a$, $\mbox{Re }a>-b$, which agree on the subset $\{a:\mbox{Re }a>0\}$.  Hence they agree on
all of $\{a:\mbox{Re }a>-b\}$, showing that $s^{-\frac{a}{2b}}U(\cp^*\cp)^{\frac{1}{2}+\frac{a}{2b}}$
is a distinguished representative of $[C_{\rho_{\eta,a}\circ\varphi}]$ when $\mbox{Re }a>-b$.  
An analogous statement holds for $[C_{\rho_{\zeta,a}\circ\sigma}]$ with $\mbox{Re }a>-c$.
The following table summarizes these conclusions.
\begin{center}
\renewcommand{\arraystretch}{1.45}
\framebox{
\begin{tabular}{c|c|c|c}
\multicolumn{4}{c}{\sc Table II} \\ \hline\hline
\multicolumn{4}{c}{Linear-fractional $\psi$ with $C_\psi$ in $C^*(C_\varphi,\mathcal{K}$)} \\ \hline
\multicolumn{1}{c}{Condition on $\psi$} &
\multicolumn{1}{|c|}{} &
\multicolumn{1}{c|}{Distinguished} &
\multicolumn{1}{c}{Matrix function}  \\[-3mm] 
\multicolumn{1}{c}{in Theorem \ref{linfrac}} &
\multicolumn{1}{|c|}{$\psi$} &
\multicolumn{1}{c|}{representative of $[C_\psi]$} &
\multicolumn{1}{c}{$\Psi(C_\psi)(t)$, $0 \leq t \leq s$} 
\\ \hline
(d) & $\begin{array}{c}\rho_{\eta,a}, \\ \mbox{Re } a > 0\end{array}$ & $s^{-\frac{a}{2b}}(C_\varphi^* C_\varphi)^{\frac{a}{2b}}$
& $\Bigg[\begin{array}{cc} \left(\frac{t}{s}\right)^{\frac{a}{2b}} & 0 \\ 0 & 0 \end{array}\Bigg]$
\\ \hline
(b) & $\begin{array}{c}\rho_{\zeta,a}, \\ \mbox{Re } a > 0\end{array}$ & $s^{-\frac{a}{2c}}(C_\varphi C_\varphi^*)^{\frac{a}{2c}}$
& $\Bigg[\begin{array}{cc} 0 & 0 \\ 0 & \left(\frac{t}{s}\right)^{\frac{a}{2c}} \end{array}\Bigg]$
\\ \hline
(a) & $\begin{array}{c}\rho_{\eta,a} \circ \varphi, \\ \mbox{Re } a > -b\end{array}$ & $s^{-\frac{a}{2b}}
U(C_\varphi^* C_\varphi)^{\frac{1}{2} + \frac{a}{2b}}$
& $\Bigg[\begin{array}{cc} 0 & \sqrt{t}\left(\frac{t}{s}\right)^{\frac{a}{2b}}  \\ 0 & 0 \end{array}\Bigg]$
\\ \hline
(c)  & $\begin{array}{c}\rho_{\zeta,a} \circ \sigma, \\ \mbox{Re } a > -c\end{array}$ & $s^{-\frac{a}{2c}-1}
U^*(\cp\cp^*)^{\frac{1}{2}+\frac{a}{2c}}$
& $\Bigg[\begin{array}{cc} 0 & 0 \\ \frac{\sqrt{t}}{s}\left(\frac{t}{s}\right)^{\frac{a}{2c}}  & 0 \end{array}\Bigg]$
\end{tabular}
}
\end{center}
Given an operator $B$ in $C^*(\cp,{\mathcal K})$, $\sigma_e(B)$ and 
$\|B\|_e$ coincide with $\sigma(\Psi(B))$ and $\|\Psi(B)\|$, respectively.  Thus, if $B$
is a finite linear combination of composition operators $C_{\psi}$ with $\psi$'s chosen from
Column 2 in Table II, one can calculate $\Psi(B)$ from
Column 4 and in principle read off $\sigma_e(B)$ and $\|B\|_e$; see Theorem 4.17 in \cite{kmm}.

It is known \cite{km} that the collection of linear-fractional composition operators $C_{\psi}$ with
$\psi$ a non-automorphism having $\|\psi\|_{\infty}=1$ is linearly independent modulo
${\mathcal K}$.  The following result shows that this remains true when ${\mathcal K}$ is
replaced by the larger subspace $C^*(C_{\varphi},{\mathcal K})$ of ${\mathcal B}(H^2)$.

\begin{thm}\label{linearcomb}
Let $\varphi$ be as in (\ref{defofph}).  Suppose that $\beta_1,\cdots,\beta_q$ are
distinct linear-fractional self-maps of $\D$ and that $a_1\cdots,a_q$ are
non-zero complex numbers.  If the linear combination $a_1C_{\beta_1}+\cdots+a_qC_{\beta_q}$
lies in $C^*(\cp,{\mathcal K})$, then so do the individual operators $C_{\beta_1},\cdots,C_{\beta_q}$. 
\end{thm}
\begin{proof}
Let us discard those $C_{\beta_i}$'s which lie in $C^*(\cp,{\mathcal K})$ and assume
for the purpose of obtaining a contradiction that there are some left over.
Relabel these as $C_{\beta_1},\cdots,C_{\beta_r}$, let $a_1,\cdots,a_r$ be the corresponding
constants, and put
$T=a_1C_{\beta_1}+\cdots+a_rC_{\beta_r}$, which lies in $C^*(\cp,{\mathcal K})$.  Here,
none of the summands are compact, so $\|\beta_i\|_{\infty}=1$ for $i=1,\cdots,r$.
Now we proceed almost as in the proof of Theorem~\ref{compacts}, with $T$ playing the role
of $C_{\psi}$.  Given $\epsilon>0$ there exists $A$ as in that proof and a complex
$c$ such that $\|T-cC_z-A\|_e<\epsilon$.  By the inequality~(\ref{lbext})
$$\|T-cC_z-A\|_e^2\geq |c|^2+\frac{1}{2\pi}\sum_{i=1}^r|a_i|^2|J(\beta_i)|$$
so that $|c|<\epsilon$, since each $|J(\beta_i)|$ must be zero, so $\beta_i$ is
a non-automorphism.  As earlier, we have $\|T-A\|_e<2\epsilon$.  

According to Corollary 5.17 in \cite{km}, the cosets $[C_{\beta_1}],\cdots,
[C_{\beta_r}]$ are linearly independent in ${\mathcal B}(H^2)/{\mathcal K}$.  
It follows that $T$ is not compact, so the matrix function
\begin{equation}
\Psi(T)=\left[\begin{array}{lr}
f_2&f_3\\f_4&f_1\end{array}\right]
\end{equation}
is not identically zero on $[0,s]$.  As in Theorem~\ref{compacts}, we focus on $f_3$ and aim for
a contradiction by assuming that $\|f_3\|_{\infty}>0$.  An appropriate choice
of $\epsilon$ again yields the inequality~(\ref{est1}), where $p$ has the same meaning as there.  
Again we write $A$ in the form~(\ref{A}), and thus 
have
\begin{eqnarray*}
4\epsilon^2>\|T-A\|_e^2&\geq&\limsup_{|z|\rightarrow 1}\left\|(T^*-A^*)\frac{k_z}{\|k_z\|}\right\|^2\\
&\geq &\lim_{\Gamma_{\zeta,D}}\left\|\left(\sum_{i=1}^r\overline{a_i}C_{\beta_i}^*-
\sum_{i=1}^m\overline{c_i}C_{\psi_i}^*\right)\frac{k_z}{\|k_z\|}\right\|^2\\
&\geq &\left\|\sum_{\zeta\in F(\psi_i)\atop D_1(\psi_i,\zeta)={\bf d_3}}\overline{c_i}k_{w_i}^+\right\|^2_{H^2_+}
\end{eqnarray*}
The rest of the argument follows that of Theorem~\ref{compacts} exactly, reaching the same contradiction.
\end{proof}

\section{Non linear-fractional maps}

In this section we explore maps $\psi$, satisfying either condition (e) or (f) of Theorem~\ref{necessity},
for which $C_{\psi}$ lies in $C^*(\cp,{\mathcal K})$; our main result shows that such maps exist.
We begin with a lemma about finite Blaschke products.
\begin{lemma}\label{finiteB}
Given $\zeta,\eta$ distinct points on $\partial \D$, and positive numbers $t_1,t_2$, there exists
a finite Blaschke product $B$ with the properties $B(\eta)=\eta,B(\zeta)=\eta,
B'(\eta)=t_1$ and $|B'(\zeta)|=t_2$.  Moreover, $B'(\zeta)=\eta\overline{\zeta}t_2$.
\end{lemma}
\begin{proof}
Clearly there is no loss of generality in taking $\eta=1$.  Initially we will
also suppose that $\zeta=-1$; this condition will be removed at the end.
A finite Blaschke product $B(z)=\prod\frac{|a_n|}{a_n}(a_n-z)/(1-\overline{a_n}z)$
will meet the conditions $B(1)=1$, $B(-1)=1$ if both
$$\prod\frac{|a_n|}{a_n}\frac{a_n-1}{1-\overline{a_n}}=1$$
and
$$\prod\frac{|a_n|}{a_n}\frac{a_n+1}{1+\overline{a_n}}=1.$$
It is easy to see that both of these conditions will be met if the zeros
of $B$ are chosen to be a collection of conjugate pairs $\{a,\overline{a}\}$.
The conditions $B'(1)=t_1, |B'(-1)|=t_2$ are satisfied if
\begin{equation}\label{der1}
\sum\frac{1-|a_n|^2}{|1-a_n|^2}=t_1
\end{equation}
and 
\begin{equation}\label{der2}
\sum\frac{1-|a_n|^2}{|1+a_n|^2}=t_2
\end{equation}
respectively (see \cite{ac}). 

Next observe that for any $t>0$, $\{z:1-|z|^2=t|1-z|^2\}$ is a circle centered
at $(t/(t+1),0)$ with radius $1/(t+1)$ and $\{z:1-|z|^2=t|1+z|^2\}$
is a circle centered at $(-t/(t+1),0)$ with radius $1/(1+t)$.  As $t\rightarrow 0$,
the centers of these circles approach $0$ and the radii tend to $1$.
Thus given $t_1, t_2$ arbitrary positive numbers we may choose $m$ a positive integer 
sufficiently large so that the circles
$\{z:1-|z|^2=\frac{t_1}{2^m}|1-z|^2\}$ and $\{z:1-|z|^2=\frac{t_2}{2^m}|1+z|^2\}$
intersect in a conjugate pair of points $a,\overline{a}$. Consider the Blaschke product $B(z)$
with a zero of order $2^{m-1}$ at $a$ and a zero of order $2^{m-1}$ at $\overline{a}$.
Since the zeros occur at conjugate pairs, $B(1)=1$ and $B(-1)=1$.
By construction 
$$(1-|a|^2)/(|1-a|^2)=(1-|\overline{a}|^2)/(|1-\overline{a}|^2)=\frac{t_1}{2^m},$$
and
$$(1-|a|^2)/(|1+a|^2)=(1-|\overline{a}|^2)/(|1+\overline{a}|^2)=\frac{t_2}{2^m},$$
so that the zeros of $B$ satisfy Equations~(\ref{der1}) and (\ref{der2}) as desired,
and $B '(1)=t_1$, $|B'(-1)|=t_2$.  

Now suppose $\zeta\in\partial \D$ is not 
equal to $-1$.  Find a parabolic automorphism
$\tau$ fixing $1$, with derivative $1$ there, and taking $\zeta$ to $-1$; a unique such $\tau$ 
exists since (purely imaginary) translations act transitively
on the boundary of the right half-plane. Then for $B$ as constructed above,
$B\circ\tau$ is a finite Blaschke product fixing $1$, with derivative $t_1$ at $1$,
sending $\zeta$ to $1$ and having derivative $|(B\circ\tau)'(\zeta)|=
|B'(-1)||\tau'(\zeta)|$; since
$|B'(-1)|$ can be arbitarity prescribed and $\tau$ depends only on
the value of $\zeta$, this means $|(B\circ\tau)'(\zeta)|$ can
be chosen to be an arbitrary positive number.  Finally observe that if $B$ is a Blaschke product with
$B(\zeta)=\eta$ and $|B'(\zeta)|=s$, then we must have $B'(\zeta)=\eta\overline{\zeta}s$,
since $\overline{\eta}\zeta B(z)$ fixes $\zeta $, and hence has positive derivative there.
\end{proof}

\begin{thm}\label{2ptfad}
Suppose that $\varphi$ is as in (\ref{defofph}).  
There exist analytic self-maps $\psi_1$ and $\psi_2$ of $\D$, satifying
conditions (e) and (f) of Theorem~\ref{necessity} respectively, such that
$C_{\psi_1}$ and $C_{\psi_2}$ lie in $C^*(C_{\varphi},{\mathcal K})$.
Moreover, $\psi_1$ and $\psi_2$
can be taken to extend continuously to the closed disk $\overline{\D}$.
\end{thm}
\begin{proof}
We shall indicate the construction for $\psi_1$ with the normalization $\eta =1$.
First consider a simply connected domain $\Omega$ in $\D$, whose boundary is a smooth
Jordan curve
which at $1$ and at $-1$ includes an arc of an internally tangent circle, such that $\overline{\Omega}\cap
\partial \D=\{-1,1\}$.  It is easy to construct a conformal map $\rho:\D$ onto
$\Omega$, extending to a homeomorphism of $\overline{D}$ onto $\overline{\Omega}$,
with the properties $\rho(1)=1, \rho(-1)=-1$ and $\rho'(1)=\frac{1}{2}$.  The map
$\rho$ will be analytic in a neighborhood of both $1$ and $-1$.  Let $\tau$ 
be the unique parabolic automorphism of $\D$ fixing $1$, having derivative $1$ there,
and mapping $\zeta$ to $-1$.  By the uniqueness statement, $\tau'(\zeta)$ is
determined by $\zeta$.  
Using Lemma~\ref{finiteB} construct a finite Blaschke product $B$ with
$B(1)=1, B(-1)=1, B'(1)=2$ and $|B'(-1)|$ chosen so that
$|B'(-1)||\rho'(-1)||\tau'(\zeta)|=s$, where $s=|\varphi'(\zeta)|$.
The map $\psi_1\equiv B\circ\rho\circ \tau$ is a self-map of $\D$, with finite angular derivative set
$F(\psi_1)=\{\zeta,1\}$ and satisfying
$\psi_1(\zeta)=1,\psi_1(1)=1, \psi_1'(1)=1 $ and
$\psi_1'(\zeta)=\overline{\zeta}s$, this last condition following from $|\psi_1'(\zeta)|=s$
and $\psi_1(\zeta)=1$.  Clearly $\psi_1$ has order of contact two at $1$ and $\zeta$.
 
Any linear-fractional map $\beta$ is uniquely determined by its second order data
vector $D_2(\beta,z_0)=(\beta(z_0),\beta'(z_0), \beta''(z_0))$ at any point
$z_0$ of analyticity.  The curvature of the curve $\psi_1(\partial \D)$ at the points
$\psi_1(\zeta)$ and $\psi_1(1)$ is determined by $D_2(\psi_1,\zeta)$ and $D_2(\psi_1,1)$, respectively.
By construction of $\psi_1$, these curvatures exceed unity.
There exists unique linear-fractional maps $\beta_1,\beta_2$ of $\D$
with the second order data vectors 
$$D_2(\beta_1,1)=D_2(\psi_1,1)=(1,1,\psi_1''(1))\mbox{ and }D_2(\beta_2,\zeta)=
D_2(\psi_1,\zeta)=(1,\overline{\zeta}s,\psi_1''(\zeta)).$$
The curvature of $\{\psi_1(e^{i\theta}):e^{i\theta}\in\partial\D\}$ matches that
of $\beta_1(\partial\D)$ at $e^{i\theta}=1$ and that of $\beta_2(\partial\D)$ at $e^{i\theta}=\zeta$.
Thus $\beta_1$ and $\beta_2$ are non-automorphism self-maps of $\D$. 
Since $\varphi'(\zeta)=\psi_1'(\zeta)=\beta_2'(\zeta)$, Theorem~\ref{linfrac} shows that 
$C_{\beta_1}$ and $C_{\beta_2}$ are in $C^*(\cp, {\mathcal K})$, and hence so is
$C_{\beta_1}+C_{\beta_2}$.  By Corollary 5.16 of \cite{km}, $C_{\psi_1}\equiv
C_{\beta_1}+C_{\beta_2} \ (\mbox{mod }{\mathcal K})$ and thus $C_{\psi_1}$
is in $C^*(\cp,{\mathcal K})$.  

For $\psi_2$, note that the Krein adjoint $\sigma$ of $\varphi$ satisfies $\sigma(\eta)=\zeta$
and $\sigma'(\eta)=1/\varphi'(\zeta)$.  We apply the first part of the proof, with $\varphi$ replaced by
$\sigma$, to find a self map $\psi_2$ of $\D$ with $F(\psi_2)=\{\zeta,\eta\}$,
$\psi_2(\eta)=\zeta$, $\psi_2'(\eta)=\sigma'(\eta)=1/\varphi'(\zeta)$, $\psi_2(\zeta)=\zeta$,
and $\psi_2'(\zeta)=1$ such that $C_{\psi_2}\in C^*(C_{\sigma},{\mathcal K})=C^*(C_{\varphi},
{\mathcal K})$ as desired.
\end{proof}

Our last theorem shows that for 
sufficiently nice $\psi$, Theorem~\ref{necessity} is almost the whole story.

\begin{thm}
Let $\psi$ be an analytic self-map of $\D$ such that $F(\psi)$ is a finite set,
$\psi$ extends analytically to a neighborhood of each point in $F(\psi)$, and for
any open set $U$ of $\partial \D$ containing $F(\psi)$, $\|\chi_{\partial \D\backslash U}\psi\|_{\infty}<1$.
If $\varphi$ is as in (\ref{defofph}), then $C_{\psi}$ lies in $C^*(\cp,{\mathcal K})$ if and only if 
\begin{itemize}
\item[(i)] one of the conditions (a)-(f) of Theorem~\ref{necessity} holds, and
\item[(ii)] the map $\psi$ has order of contact two at each point of $F(\psi)$.
\end{itemize}
\end{thm}
\begin{proof}
Suppose $\psi$ is as described and $C_{\psi}$ is in $C^*(\cp,{\mathcal K})$. By Theorem~\ref{necessity},
one of (a)-(f) holds.  Suppose that $\gamma$ is in $F(\psi)$ and $\psi$ has order of contact
exceeding two at $\gamma$. Following Theorem~\ref{compacts}, we let $\epsilon>0$ and
find a finite linear combination $A$ of composition operators whose self-maps are chosen from
the lists (\ref{lists}) such that $\|C_{\psi}-A\|<\epsilon$.  At the same time we apply
the inequality ~(\ref{lb2}) to the linear combination $C_{\psi}-A$ at the point $\alpha=\gamma$.
The maps in the lists (\ref{lists}), being linear-fractional non-automorphisms, all have
order of contact two at the unique points in their angular derivative sets.  Taking the left
side of (\ref{lb2}) to be $\|C_{\psi}-A\|_e^2$ and $k=3$
on the right side, the sum on the right-hand side has only one term, namely $1/|\psi'(\gamma)|$,
giving
$$\epsilon^2>\|C_{\psi}-A\|_e^2\geq\frac{1}{|\psi'(\gamma)|},$$
a contradiction.  Thus $\psi$ must have order of contact two at $\gamma$.

Conversely, suppose $\psi$ satisfies (i) and (ii).   If $F(\psi)=\{\zeta\}$, let $\beta$ be the 
unique linear-fractional map with $D_2(\beta,\zeta)=D_2(\psi,\zeta)$.  Since $\psi$ has
order of contact two at $\zeta$, the curvature of $\{\psi(e^{i\theta}):e^{i\theta}\in\partial \D\}$
at $e^{i\theta}=\zeta$ exceeds unity, so that $\beta(\partial\D)$, having the same curvature, is
internally tangent to $\partial\D$ at $\psi(\zeta)$ and bounds a proper subdisk of $\D$. Thus
$\beta$ is a non-automorphism of $\D$ satisfying (a) or (b) of Theorem~\ref{linfrac}, and $C_{\beta}$
lies in $C^*(\cp,{\mathcal K})$.  The same argument covers the case $F(\psi)=\{\eta\}$.

If $F(\psi)=\{\zeta,\eta\}$ we proceed as in the proof of Theorem~\ref{2ptfad} to produce linear-fractional
non-automorphisms $\beta_1$ and $\beta_2$ of $\D$ with
$C^*(\cp,{\mathcal K})$ containing $C_{\beta_1}$ and $C_{\beta_2}$ and 
$C_{\psi}\equiv C_{\beta_1}+C_{\beta_2}\ (\mbox{mod }{\mathcal K})$. 
\end{proof}


\begin{thebibliography}{9}






\bibitem{blns}
P. Bourdon, D. Levi, S. Narayan and J. Shapiro,
Which linear fractional composition operators are essentially normal?,
J. Math. Anal. Appl. 280 (2003), 30--53.


\bibitem{ac}
P. Ahern and D. Clark, On inner functions with $H^p$-derivative,
Michigan Math. J. 21 (1974), 115--127.

\bibitem{br}
E. Basor and D. Retsek, Extremal non-compactness of composition operators
with linear fractional symbol, J. Math. Anal. Appl. 322 (2006), 749--763.

\bibitem{b} 
E. Berkson, Composition operators isolated in the uniform operator
topology, Proc. Amer. Math. Soc. 81 (1981), 230--232.

\bibitem{cocot}
J. Conway, \textit{A Course in Operator Theory},
Amer. Math. Soc., Providence, 2000.
 
\bibitem{cowen}
C. Cowen, Composition operators on $H^2$, J. Operator Theory 9 (1983),
77--106.

\bibitem{cow2}
C. Cowen, Linear-fractional composition operators on $H^2$, Integral
Equations Operator Theory 11 (1988), 151--160.

\bibitem{cmbook}
C. Cowen and B. MacCluer,
\textit{Composition Operators on Spaces of Analytic Functions},
CRC Press, Boca Raton, 1995.


\bibitem{dash}
A. Dash, Joint essential spectra, Pacific J. Math., 64 (1976), 119--128.


\bibitem{g}
J. Guyker, On reducing subspaces of composition operators,
Acta Sci. Math. (Szeged) 53 (1989), 369--376.

\bibitem{HP}
E. Hille and R. Phillips, Functional Analysis and Semi-Groups, 
Amer. Math. Soc., Providence, 1957.

\bibitem{j}
M. Jury, $C^*$-algebras generated by groups of composition operators, 
Indiana Univ. Math. J. 56 (2007), 3171--3192.

\bibitem{j2}
M. Jury, The Fredholm index for elements of Toeplitz-composition $C^*$-algebras,
Integral Equations Operator Theory, 58 (2007), 341--362.


\bibitem{km}
T. Kriete and J. Moorhouse, Linear relations in the Calkin algebra
for composition operators, Trans. Amer. Math. Soc. 359 (2007), 2915--2944.

\bibitem{kmm}
T. Kriete, B. MacCluer, and J. Moorhouse, Toeplitz-composition $C^*$-algebras,
J. Operator Theory 58 (2007), 135--156.

\bibitem{kmm3}
T. Kriete, B. MacCluer, and J. Moorhouse, Spectral theory for algebraic
combinations of Toeplitz and composition operators, preprint.

\bibitem{ms}
B. MacCluer and J. Shapiro, Angular derivatives and compact composition
operators on the Hardy and Bergman spaces, Canadian J. Math. 38 (1986), 878--906.

\bibitem{montes}
A. Montes-Rodr\'iguez, M. Ponce-Escudero, and S. Shkarin,
Invariant subspaces of parabolic self-maps in the Hardy space,
preprint.


\bibitem{ss}
J. Shapiro and C. Sundberg, Isolation amongst the composition operators,
Pacific J. Math. 145 (1990), 117--152.

\bibitem{st}
J. Shapiro and P. Taylor, Compact, nuclear, and Hilbert-Schmidt composition
operators on $H^2$, Indiana Univ. Math. J. 23 (1973), 471--496.


\bibitem{sbook}
J. Shapiro, 
\textit{Composition Operators and Classical Function Theory}, 
Springer-Verlag, New York, 1993.





\end{thebibliography}
\end{document}